\documentclass[11pt]{amsart}
\usepackage{CJK}
\usepackage{amssymb}
\usepackage{amsmath}
\usepackage{amsthm}

\usepackage{graphicx}
\usepackage{pst-plot}

\newtheorem{Def}{Definition}
\newtheorem{Prop}[Def]{Proposition}
\newtheorem{Thm}[Def]{Theorem}
\newtheorem{Lem}[Def]{Lemma}
\newtheorem{Cor}[Def]{Corollary}

\newtheorem*{Thm0}{Theorem}

\theoremstyle{definition}
\newtheorem{Ex}[Def]{Example}

\begin{document}

\title{On the first singularity for the Upsilon invariant of algebraic knots}

\author{Shida Wang}

\address{Department of Mathematics, Indiana University, Bloomington, IN 47405}
\email{wang217@indiana.edu}

\maketitle

\begin{abstract}

We show that the location of the first singularity of the Upsilon function of an algebraic knot is determined by the first term of its Puiseux characteristic sequence.
In many cases this gives better bounds than the tau invariant on the genus of a cobordism between algebraic knots.

\end{abstract}

\section{Introduction}

The slice genus, or 4--genus, $g_4$ of a knot $K$ is the minimal genus of a smoothly embedded oriented surface in $B^4$ whose boundary is $K$.
A cobordism between two knots, $K_0,K_1\subset S^3$, is a smoothly embedded oriented surface in $S^3\times[0,1]$ with boundary $K_0\times\{0\}\cup-K_1\times\{1\}$.
It is easy to see that the genus of any cobordism between $K_0$ and $K_1$ is bounded below by $|g_4(K_1)-g_4(K_0)|$.
If a cobordism of genus $|g_4(K_1)-g_4(K_0)|$ exists, it is called \emph{minimal}.
A natural question asks, for given knots, whether a minimal cobordism exist.
See \cite{scissor,Gordian,minimal,immersed} for related work.

For algebraic knots, if the cobordism is algebraic,
that is, if the cobordism is the transverse intersection of a smooth algebraic curve with the compact domain enclosed by two spheres of different radii in $\mathbb{C}^2$,
then a corollary of the Thom conjecture asserts it is minimal.
Even if an algebraic cobordism does not exist, there might be a smooth one which is minimal.
In this note we will give new families of examples that do not have smooth minimal cobordisms at all.

There is an obvious relationship between a minimal cobordism and a minimal unknotting sequence:
if $K_1$ is a knot with the same 4--genus and unknotting number and $K_0$ appears in some minimal unknotting sequence for $K_1$, then there is a minimal cobordism between $K_0$ and $K_1$.
In \cite{main}, Borodzik and Livingston proved the following result by semicontinuity of semigroups.

\begin{Thm0} \emph{(\cite[Theorem 2.17]{main})} If the torus knot $T_{a,b}$ with $0<a<b$ appears in some minimal unknotting sequence for the torus knot $T_{c,d}$ with $0<c<d$, then $a\leqslant c$.\end{Thm0}

This theorem originates from the nonexistence of a positively minimal self-intersecting concordance, but it does not imply the nonexistence of a minimal cobordism, which we will prove as an easy corollary in this note.

\begin{Thm0}If there is a minimal cobordism between the torus knots $T_{a,b}$ with $0<a<b$ and $T_{c,d}$, and $g_4(T_{a,b})\leqslant g_4(T_{c,d})$, then $a\leqslant c$.\end{Thm0}
\vspace{-1em}

Moreover, we will express this generally for algebraic knots (see Corollary \ref{cobordism}).

A numerical invariant $\tau$ was introduced in \cite{tau}.
It gives a homomorphism from the smooth concordance group to $\mathbb{Z}$ and satisfies $|\tau|\leqslant g_4$.
For any positive torus knot $K$, $\tau(K)=g(K)=g_4(K)$, where $g$ is the 3--genus.
More generally, this holds for L-space knots (\cite[Theorem 1.2]{Lknot} combined with \cite[Theorem 1.2]{HFKgenus}), which include algebraic knots (cf. \cite[Theorem 1.10]{cableL}).
The $\tau$ invariant is successful in showing $g=g_4$ for algebraic knots,
but cannot obstruct a cobordism of genus $|g_4(K_1)-g_4(K_0)|$ for algebraic knots $K_0$ and $K_1$, because $|\tau(K_1\#-K_0)|=|\tau(K_1)-\tau(K_0)|=|g_4(K_1)-g_4(K_0)|$.

\vspace{-0.5em}

Although $\tau$ cannot be used to prove the above theorem, one may wonder if the Tristram-Levine signature function can.
We will give a counterexample showing it cannot (see the remark after Corollary \ref{cobordism}).

The invariant we will use to obstruct the minimal cobordism is called the Upsilon function, which generalizes $\tau$.
In \cite{upsilon} this function, $\Upsilon_K(t)\colon[0,2]\rightarrow\mathbb{R}$, is associated to any knot $K\subset S^3$.
It is derived from the knot Heegaard Floer homology and it is shown to be a concordance invariant.
It satisfies the symmetry property $\Upsilon_K(2-t)=\Upsilon_K(t)$.
The important property we will use is the following.

\begin{Prop}\label{genusbound} \emph{(\cite[Theorem 1.11]{upsilon})} The invariants $\Upsilon_K(t)$ bound the 4--genus of $K$;
more precisely, for $0\leqslant t\leqslant1$, $|\Upsilon_K(t)|\leqslant t\cdot g_4(K)$.\end{Prop}

\vspace{-0.5em}

Note that $\Upsilon$ is a piecewise linear function and its slope near $t=0$ is $-\tau$.
It then follows that the Upsilon function near $t=0$ does not obstruct minimal cobordisms.
Since Ozsv\'{a}th, Stipsicz and Szab\'{o} show $\Upsilon$ is fully computable for L-space knots in terms of the Alexander polynomial (\cite[Theorem 6.2]{upsilon}),
we will use such a computation to determine the value of $\Upsilon$ to the right of its first singularity
so that a stronger obstruction than $\tau$ is obtained.

\vspace{-0.5em}

This note is organized as follows. Section 2 gives a brief account of the theory of algebraic knots and iterated torus knots.
Section 3 computes the first singularity of the Upsilon function of an algebraic knot in terms of its classical invariants.
Section 4 uses the results of Section 3 to obstruct minimal cobordisms and some examples are given.

\vspace{-0.5em}

\emph{Acknowledgments.} The author wishes to express sincere thanks to Professor Charles Livingston for proposing this study and carefully reading a draft of this paper.
Thanks also to Professor Maciej Borodzik for valuable comments.
The author recently learned that similar results to those presented here have been obtained by Peter Feller independently.

\section{Algebraic knots as iterated torus knots}

We refer to \cite[\S10]{Puiseux1,Puiseux2,Puiseux3} for a general introduction to singular points of complex curves.

By an isolated plane curve singular point we mean a pair $(C,z)$, where $C$ is a complex curve in $\mathbb{C}^2$ and $z$ is a point on $C$,
such that $C$ is smooth at all points sufficiently close to $z$, with the exception of $z$ itself.
For a sufficiently small $r>0$, $C$ intersects the ball $B(z,r)\subset\mathbb{C}^2$ transversally along a link $L$ and $C\cap B(z,r)$ is topologically a cone on this link.
The link $L$ will be called the link of the singular point.
A singular point is called cuspidal or unibranched if $L$ is connected, that is, if $L$ is a knot.
A knot that arises in this way is called an \emph{algebraic knot}.

Any cuspidal singular point is topologically equivalent to a singular point which is locally parametrized by $x=t^{q_{_0}},y=t^{q_{_1}}+\cdots+t^{q_{_n}}$,
where $q_{_0}<q_{_1}<\cdots<q_{_n}$ are positive integers.
There is a unique such representation if we further assume $\gcd(q_{_0},\cdots,q_{_i})$, denoted by $D_i$, does not divide $q_{_{i+1}}$ and $D_n=1$.
Such a sequence $(q_{_0};q_{_1},\cdots,q_{_n})$ is called the \emph{Puiseux characteristic sequence} of the singular point.
Sometimes the first term $q_{_0}$ of this sequence is called the \emph{multiplicity} of the singular point.

For any singular point $(C,z)$ one can associate a semigroup $S\subset\mathbb{Z}_{\geqslant0}$.
Let $\varphi(t)=(x(t),y(t))$ be a local parametrization of $C\cap B(z,r)$ near $z$;
that is, $\varphi$ is analytic and gives a homeomorphism from an open neighborhood of $0\in\mathbb{C}$ to $C\cap B(z,r)$ such that $\varphi(0)=z$.
Then $\varphi$ induces a map $\varphi^*\colon\mathbb{C}[[x,y]]\rightarrow\mathbb{C}[[t]]$ by $f(x,y)\mapsto f(x(t)-z_1,y(t)-z_2)$, where $(z_1,z_2)=z\in\mathbb{C}^2$.
The map $\mathrm{ord}\colon\mathbb{C}[[t]]\rightarrow\mathbb{Z}_{\geqslant0}\cup\{\infty\}$ maps a power series in one variable to its order at 0 (the order of the zero function is defined to be $\infty$).
The set $S$ of nonnegative integers belonging to the image of $\mathbb{C}[[x,y]]$ under the composition $\mathrm{ord}\circ\varphi^*$ is easily seen to be a semigroup,
and it is defined to be the \emph{semigroup} of the cuspidal singular point \cite[\S 4.3]{Puiseux1}.
Sometimes we abbreviate this terminology and call it the semigroup of the algebraic knot.

The link of a cuspidal singular point with Puiseux characteristic sequence\linebreak $(q_{_0};q_{_1},\cdots,q_{_n})$
is an $(n-1)$-fold iterate of a torus knot $T_{\frac{q_{_0}}{D_1},\frac{q_{_1}}{D_1}}$,
whose parameters are determined by $(q_{_0};q_{_1},\cdots,q_{_n})$. Also, the semigroup $S$ is determined by $(q_{_0};q_{_1},\cdots,q_{_n})$. The details are as follows.

Define $s_0=q_{_0}$ and $s_i=\frac{D_0q_{_1}+D_1(q_{_2}-q_{_1})+\cdots+D_{i-1}(q_{_i}-q_{_{i-1}})}{D_{i-1}}$ (cf. \cite[(4.4)(4.5)]{Puiseux1}).

\begin{Prop}\label{translate} Let $(C,z)$ be a cuspidal singular point with Puiseux characteristic sequence $(q_{_0};q_{_1},\cdots,q_{_n})$ and $K$ be its link.

(i) \emph{(\cite[p.83 Theorem 4.3.5(i)]{Puiseux1})} The semigroup $S$ is generated by $s_0,s_1,\cdots,s_n$.

(ii) \emph{(\cite[p.121 or Theorem 9.8.1]{Puiseux1})} The knot $K$ is an iterated torus knot
$$\Big(\Big(\mbox{\LARGE T}_{\frac{D_0}{D_1},\frac{s_1}{D_1}}\Big)\,_{\frac{D_1}{D_2},\frac{s_2}{D_2}}\Big)_{\cdots\frac{D_{n-1}}{D_n},\frac{s_n}{D_n}}\raisebox{-1em}{.}$$

(iii) \emph{(\cite[Theorem 6.3.2]{Puiseux1} or \cite[Remark 10.10]{Puiseux3})} The genus of $K$ is $\frac{1}{2}\mu$, where $\mu$ is the Milnor number and equals $\sum_{i=1}^n(q_{_i}-1)(D_{i-1}-D_i)$.

\end{Prop}

\begin{Ex} (1) Let $K$ be an algebraic knot with Puiseux characteristic sequence $(4;6,7)$.
Then $D_0=4,D_1=2,D_2=1$.
So $s_0=4,s_1=6,s_2=13$.
Thus $K=(T_{2,3})_{2,13}$.
It has genus 8.

(2) Let $K$ be an algebraic knot with Puiseux characteristic sequence $(12;18,22,25)$.
Then $D_0=12,D_1=6,D_2=2,D_3=1$.
So $s_0=12,s_1=18,s_2=40,s_3=123$.
Thus $K=((T_{2,3})_{3,20})_{2,123}$.
It has genus 105.
\end{Ex}

Finally we remind the readers that $g=g_4$ for algebraic knots \cite{gauge}. A knot Floer approach can be found in \cite{tau}.

\section{The Upsilon function of algebraic knots}

We refer to \cite{upsilon} for the definition of the Upsilon invariant.
For our purpose, we only need to know that $\Upsilon$ is a piecewise linear function and its slope near $t=0$ is $-\tau$.
Also, $|\tau|=g=g_4$ for algebraic knots \cite{tau}, as mentioned in Section 1.

In \cite[Theorem 6.2]{upsilon}, the Upsilon invariant of L-space knots is computed in terms of the Alexander polynomial.
Alternatively, it can be expressed in terms of semigroups for algebraic knots as follows.

\begin{Prop}\label{formula} \emph{(\cite[Proposition 3.4]{main})} Let $K$ be an algebraic knot with genus $g$ and $S$ be the corresponding semigroup.
Then for any $t\in[0,1]$ we have $$\Upsilon_K(t)=\max\limits_{m\in\{0,\cdots,2g\}}\{-2\#(S\cap[0,m))-t(g-m)\}.$$\end{Prop}

\emph{Remark.} Recently, \cite{Legendre} shows that the Upsilon function of an algebraic knot is a Legendre transform of a function determined by algebraic invariants of the knot,
related to $d$-invariants of surgery on $K$.

We will prove the following

\begin{Thm}\label{singularity}Let $K$ be an algebraic knot with genus $g$ and the first term of its Puiseux characteristic sequence $a$.
Then $\Upsilon_K(t)=-gt$ for $t\in[0,\frac{2}{a}]$ and $\Upsilon_K(t)>-gt$ for $t>\frac{2}{a}$.\end{Thm}

Thus the shape of Upsilon function for an algebraic knot looks like the following graph.

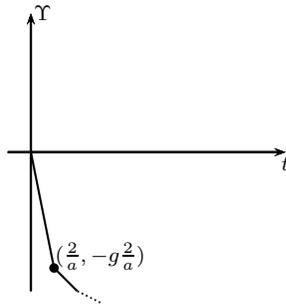
\vspace{-1em}

\psset{unit=4em}

\rput(3,-1){\psaxes[ticks=none,labels=none]{->}(0,0)(-0.2,-1.2)(2.2,1.2)
\rput(2.2,-0.1){\scriptsize$t$}\rput(0.1,1.2){\scriptsize$\Upsilon$}
\psdots(0.2,-1)\rput(0.6,-0.9){\scriptsize$(\frac{2}{a},-g\frac{2}{a})$}
\psplot[plotstyle=curve]{0}{0.2}{x -5 mul}
\psplot[plotstyle=curve]{0.2}{0.4}{x -1 mul 0.8 sub}
\psline[linestyle=dotted,dotsep=1pt](0.4,-1.2)(0.6,-1.3)
}

\vspace{8em}

\ \

First we prove a lemma.

\begin{Lem}\label{semigroup}Let $K$, $a$ and $S$ be as above. Then $a\cdot\#(S\cap[0,m))\geqslant m$.\end{Lem}

\textbf{Proof.} By Proposition \ref{translate}(i), $a$ is a generator of $S$.
Since $a$ is a positive integer, there is a nonnegative integer $n$ such that $na<m\leqslant(n+1)a$.
Clearly $0,a,2a,\cdots,na\in S\cap[0,m)$.
Thus $\#(S\cap[0,m))\geqslant n+1$ and therefore $a\cdot\#(S\cap[0,m))\geqslant a(n+1)\geqslant m$.\hfill$\Box$

\textbf{Proof of Theorem \ref{singularity}.} Taking $m=0$, we have the linear function $-2\#(S\cap[0,m))-t(g-m)=-gt$. So $\Upsilon_K(t)\geqslant-gt$.
On the other hand, Lemma \ref{semigroup} implies $-2\#(S\cap[0,m))-t(g-m)\leqslant-2\frac{m}{a}-tg+tm\leqslant-gt$ if $t\leqslant\frac{2}{a}$.
Hence $\Upsilon_K(t)\geqslant-gt$ for $t\in[0,\frac{2}{a}]$.

To prove the second part of the theorem, note that by \cite[p.86 line 19]{Puiseux1} $a$ is the least nonzero element of $S$.
Thus, $\#(S\cap[0,m))=1$ for $1\leqslant m\leqslant a$.
In particular, $-2\#(S\cap[0,m))-t(g-m)=-2-t(g-a)$ when $m=a$.
Here we can take $m=a$ in Proposition \ref{formula} because $a\leqslant2g$ follows easily from Proposition \ref{translate}(iii).
Evaluating the linear function $-2-t(g-a)$ on the interval $(\frac{2}{a},1]$ yields $\Upsilon_K(t)>-gt$.\hfill$\Box$

Let $T_{a,b}$ be a torus knot with $a<b$. Then it is an algebraic knot with Puiseux characteristic sequence $(a;b)$.

\begin{Cor}For any torus knot $T_{a,b}$ with $a<b$,
we have $\Upsilon_{T_{a,b}}(t)=-g(T_{a,b})t$ for $t\in[0,\frac{2}{a}]$ and $\Upsilon_{T_{a,b}}(t)>-g(T_{a,b})t$ for $t>\frac{2}{a}$.\end{Cor}

\section{Applications to minimal cobordisms between algebraic knots}

Using Theorem \ref{singularity}, we easily get

\begin{Cor}\label{cobordism}Let $K_0$ and $K_1$ be two algebraic knot with the first term of their Puiseux characteristic sequence $a$ and $c$, respectively.
If $a>c$ and $g(K_0)\leqslant g(K_1)$, then $g_4(K_1\#-K_0)>g(K_1)-g(K_0)$; that is, there is no minimal cobordism between $K_0$ and $K_1$.\end{Cor}

\textbf{Proof.} Evaluating at $t=\frac{2}{c}$ (which is greater than $\frac{2}{a}$),
$\Upsilon_{K_0}(\frac{2}{c})>-g(K_0)\frac{2}{c}$ and $\Upsilon_{K_1}(\frac{2}{c})=-g(K_1)\frac{2}{c}$.
Thus $\Upsilon_{K_1\#-K_0}(\frac{2}{c})=\Upsilon_{K_1}(\frac{2}{c})-\Upsilon_{K_0}(\frac{2}{c})<-g(K_1)\frac{2}{c}+g(K_0)\frac{2}{c}$.
Hence $\frac{2}{c}g_4(K_1\#-K_0)\geqslant|\Upsilon(K_1\#-K_0)(\frac{2}{c})|>\frac{2}{c}(g(K_1)-g(K_0))$.\hfill$\Box$

\emph{Remark.} This result cannot be obtained by using the Tristram-Levine signature function as an obstruction.
For example, Corollary \ref{cobordism} implies that there is no minimal cobordism between the torus knots $T_{4,5}$ and $T_{3,10}$.
However, using a combinatorial formula \cite[Proposition 1]{signature} for the Tristram-Levine signature functions of torus knots (also see, for example, \cite[Lemma 15]{Gordian} or \cite[(1)]{scissor}),
one computes that the absolute value of the signature function of $T_{4,5}\#-T_{3,10}$ is bounded above by 6.
The strongest consequence of this is that the minimal genus of cobordisms between $T_{4,5}$ and $T_{3,10}$ is 3, which equals $g_4(T_{3,10})-g_4(T_{4,5})$.

Since algebraic cobordisms are minimal, we have

\begin{Cor}Let $K_0$ and $K_1$ be two algebraic knot with the first term of their Puiseux characteristic sequence $a$ and $c$, respectively.
If $g_4(K_0)\leqslant g_4(K_1)$ and there is an algebraic cobordism between $K_0$ and $K_1$, then $a\leqslant c$.\end{Cor}

Due to the relationship between a minimal cobordism and a minimal unknotting sequence mentioned in the introduction, we have

\begin{Cor}Let $K_0$ and $K_1$ be two algebraic knot with the first term of their Puiseux characteristic sequence $a$ and $c$, respectively.
If $K_0$ appears in some minimal unknotting sequence of $K_1$, then $a\leqslant c$.\end{Cor}

\textbf{Proof.} For algebraic knots, the unknotting number $u$ is bounded above by $\frac{\mu}{2}=g$. See \cite[Remark 10.10]{Puiseux3} and \cite[\S 4]{unknotting}.
The trace of the minimal unknotting sequence gives a genus $u$ smooth surface in $B^4$ bounded by the knot, so $u\geqslant g_4$.
By \cite{gauge}, $g_4=g$ for algebraic knots. Hence $u=g_4$ for algebraic knots.

Let $\Sigma_1$ be a genus $u(K_1)$ smooth surface in $B^4$ given by the minimal unknotting sequence of $K_1$ in the hypothesis.
Then there is a subsurface $\Sigma_0\subset\Sigma_1$ bounded by $K_0$.
Clearly $g(\Sigma_0)\geqslant g_4(K_0)$.
The surface $\Sigma_1\backslash\Sigma_0$ gives a cobordism of genus less than or equal to $g_4(K_1)-g_4(K_0)$.
This cannot happen if $a\leqslant c$.\hfill$\Box$

\begin{Cor} \emph{(\cite[Theorem 2.17]{main})} If the torus knot $T_{a,b}$ with $0<a<b$ appears in some minimal unknotting sequence
for the torus knot $T_{c,d}$ with $0<c<d$, then $a\leqslant c$.\end{Cor}

\emph{Remark.} Although related, minimal unknotting sequences for torus knots do not necessarily yield minimal cobordisms.
By \cite[Theorem 3 and Proposition 23]{Gordian}, there exists a minimal cobordism between $T_{2,3c}$ and $T_{3,2c}$, but $T_{2,3c}$ is not in any minimal unknotting sequence of $T_{3,2c}$, for any $c>1$.

\begin{Ex}Let $a,b,c$ satisfy $2a<b<\frac{c}{3}$, $\gcd(a,b)=\gcd(2a,b)=\gcd(a,3b)=d>1$ and $\gcd(a,b,c)=1$.
Consider algebraic knots $K_0$ and $K_1$ with Puiseux characteristic sequences $(2a,b,c)$ and $(a,3b,c)$, respectively.
It is easy to verify that $g_4(K_0)\leqslant g_4(K_1)$.
Thus there is no algebraic cobordism between $K_0$ and $K_1$.
Also $K_0$ cannot appear in any minimal unknotting sequence of $K_1$.

For example, if we take $a=4,b=10,c=31$, then $K_0$ and $K_1$ have Puiseux characteristic sequences $(8,10,31)$ and $(4,30,31)$, respectively.
This means $K_0=(T_{4,5})_{2,61}$ and $K_1=(T_{2,15})_{2,61}$.
Hence $(T_{4,5})_{2,61}$ cannot appear in any minimal unknotting sequence of $(T_{2,15})_{2,61}$.
\end{Ex}


\begin{thebibliography}{Hed10}
\bibitem[Bad12]{scissor}S. Baader, \emph{Scissor equivalence for torus links}, Bull. London Math. Soc. 44 (2012) 1068--1078.
\bibitem[BM84]{unknotting}M. Boileau and C. Weber, \emph{le probl\`{e}me de J. Milnor sur le nombre gordien des n{\oe}uds alg\'{e}briques}, l'Enseignement math\'{e}matique, t.30 (1984) 173--222.
\bibitem[BH15]{Legendre}M. Borodzik and M. Hedden, \emph{The $\Upsilon$ function of $L$-space knots is a Legendre transform}, arXiv:1505.06672v1.
\bibitem[BL13]{main}M. Borodzik and C. Livingston, \emph{Semigroups, d-invariants and deformations of cuspidal singular points of plane curves}, arXiv:1305.2868v2.
\bibitem[EN85]{Puiseux2}D. Eisenbud and W. Neumann, \emph{Three-dimensional link theory and invariants of plane curve singularities}, Annals of Mathematics Studies 110, Princeton Universityty Press, Princeton, 1985.
\bibitem[Fel14]{Gordian}P. Feller, \emph{Gordian adjacency for torus knots}, Algebraic \& Geometric Topology 14 (2014) 769--793.
\bibitem[Fel15]{minimal}P. Feller, \emph{Minimal cobordisms between torus knots}, arXiv:1501.00483v2.
\bibitem[Hed10]{cableL}M. Hedden, \emph{On Knot Floer Homology and Cabling: 2}, International Mathematics Research Notices, Vol. 2009, No. 12, pp. 2248–-2274.
\bibitem[KM93]{gauge}P. B. Kronheimer and T. S. Mrowka, \emph{Gauge theory for embedded surfaces, I}, Topology 32, No. 4 (1993) 773--826.
\bibitem[Lit79]{signature}R. A. Litherland, \emph{Signatures of iterated torus knots}, Topology of Low-Dimensional Manifolds (Proceedings of the Second Sussex Conference, 1977) pp. 71--84,
Lecture Notes in Mathematics, Volume 722, Springer, Berlin, 1979.
\bibitem[Mil68]{Puiseux3}J. Milnor, \emph{Singular points of complex hypersurfaces}, Annals of Mathematics Studies 61, Princeton University Press and the University of Tokyo Press, Princeton, New Jersey, 1968.
\bibitem[OSt13]{immersed}B. Owens and S. Strle, \emph{Immersed disks, slicing numbers and concordance unknotting numbers}, arXiv:1311.6702v2.
\bibitem[OSS14]{upsilon}P. Ozsv\'{a}th, A. Stipsicz and Z. Szab\'{o}, \emph{Concordance homomorphisms from knot Floer homology}, arXiv:1407.1795v2.
\bibitem[OSz03]{tau}P. Ozsv\'{a}th and Z. Szab\'{o}, \emph{Knot Floer homology and the four-ball genus}, Geometric \& Topology 7 (2003) 615--639.
\bibitem[OSz04]{HFKgenus}P. Ozsv\'{a}th and Z. Szab\'{o}, \emph{Holomorphic disks and genus bounds}, Geometric \& Topology 8 (2004) 311--334.
\bibitem[OSz05]{Lknot}P. Ozsv\'{a}th and Z. Szab\'{o}, \emph{On knot Floer homology and lens space surgeries}, Topology 44 (2005) 1281--1300.
\bibitem[Wal04]{Puiseux1}C. T. C. Wall, \emph{Singular Points of Plane Curves}, London Mathematical Society Student Texts 63, Cambridge University Press, Cambridge, 2004.
\end{thebibliography}
\end{document}